\documentclass[preprint,12pt]{elsarticle}

\usepackage{amsthm,amssymb,multirow,mathtools}

\usepackage[all]{xy}
\usepackage[bookmarks,colorlinks]{hyperref}

\usepackage[linesnumbered,ruled,vlined]{algorithm2e}

\numberwithin{equation}{section}

\newtheorem{theorem}{Theorem}

\theoremstyle{definition}
\newtheorem{definition}{Definition}[section]
\newtheorem{lemma}[definition]{Lemma}
\newtheorem{proposition}[definition]{Proposition}
\newtheorem{corollary}[definition]{Corollary}
\newtheorem{example}[definition]{Example}
\newtheorem{remark}[definition]{Remark}

\newcommand{\A}{\mathcal{A}}
\newcommand{\ad}{\mathrm{ad}}
\newcommand{\aff}{\mathfrak{aff}}
\newcommand{\Au}{\mathrm{Aut}}
\newcommand{\B}{\mathcal{B}}
\newcommand{\C}{\mathbb C}
\newcommand{\Cc}{\mathcal C}
\newcommand{\co}{\coloneqq}
\newcommand{\Der}{\mathrm{Der}}
\newcommand{\GL}{\mathrm{GL}}
\newcommand{\id}{\mathrm{id}}
\newcommand{\im}{\mathrm{im}\,}
\newcommand{\K}{\mathbb K}
\newcommand{\N}{\mathbb N}
\newcommand{\M}{\mathrm{Mat}}
\newcommand{\R}{\mathbb R}
\newcommand{\rc}{\mathcal R}
\newcommand{\s}{\mathrm{span}}
\newcommand{\Sc}{\mathcal S}
\newcommand{\Sf}{\mathfrak s}
\newcommand{\Z}{\mathcal Z}

\begin{document}

\begin{frontmatter}

\title{Complete and cocomplete Lie algebras with injective and projective properties}

\author[1,2]{Vu A. Le}
\ead{vula@uel.edu.vn}
\address[1]{University of Economics and Law, Ho Chi Minh City, Vietnam}
\address[2]{Vietnam National University, Ho Chi Minh City, Vietnam}

\author[3]{Hoa Q. Duong}
\ead{hoa.duongquang@ufm.edu.vn}
\address[3]{University of Finance - Marketing, Ho Chi Minh City, Vietnam}

\author[4]{Tuan A. Nguyen\corref{cor1}}
\ead{tuannguyenanh@hcmue.edu.vn}
\address[4]{Ho Chi Minh City University of Education, Ho Chi Minh City, Vietnam}
\cortext[cor1]{Corresponding author: tuannguyenanh@hcmue.edu.vn}

\begin{abstract}
Motivated by the classical correspondence between short exact sequences and splitting properties in module theory, this paper investigates projective and injective analogues in the category of Lie algebras. We first establish that no Lie algebra can serve as a projective or injective object with respect to arbitrary extensions, thus clarifying the natural limitation of this analogy. To recover meaningful dual behaviors, we introduce two new notions: \emph{cocentral extensions} and \emph{cocomplete Lie algebras}, viewed as the natural dual counterparts of central extensions and complete Lie algebras. We prove that solvable complete Lie algebras exhibit an injective-like property, while cocomplete Lie algebras $\Cc$ satisfy the cohomological condition $H^2(\Cc, \K) = 0$, acting as projective-like objects. Moreover, we obtain a full classification of almost abelian cocomplete Lie algebras. These results establish a ``duality framework’’ for \emph{completeness} and \emph{cocompleteness} in Lie algebra extensions, linking structural, categorical, and cohomological aspects.
\end{abstract}


\end{frontmatter}

\tableofcontents

\section{Introduction}\label{sec1}

Throughout this paper, unless otherwise stated, $n$ denotes a positive integer, all considered Lie algebras and vector spaces are finite-dimensional over a field $\K$ of characteristic 0. For convenience, an $n$-dimensional $\K$-vector space is always presented by $\K^n$ and $V$ denotes some finite-dimensional $\K$-vector space. The notations $\GL _n(\K)$ and $\mathfrak{gl}_n(\K) \equiv \M_n(\K)$ 
denote the Lie group of invertible $(n \times n)$-matrices and the Lie algebra of $(n\times n)$-matrices with entries in $\K$, respectively. We also denote by $\A \oplus \B$ the direct sum of Lie algebras $\A$ and $\B$, while $\A \oplus _\rho \B$ stands for the semidirect sum of $\A$ and $\B$ via an action $\rho \colon \B \to \Der(\A)$ of $\B$ on $\A$. Besides, $\s \{x_1, \ldots, x_n\}$ indicates an $n$-dimensional vector space with basis $(x_1, \ldots, x_n)$ and $U\dotplus V$ is the direct sum of vector spaces $U$ and $V$.

This paper is concerned with short exact sequences of Lie algebras 
\begin{equation}\label{LE0}
    \xymatrix
	{
 		0 \ar[r] & \A \ar[r]^\iota & \B \ar[r]^\pi & \Cc \ar[r] & 0
 	}\tag{LE}
\end{equation}
which is called \emph{Lie algebra extensions} (for short, \emph{extensions}) of $\Cc$ by $\A$. We are interested in structural conditions on Lie algebras that force such extensions to \emph{split trivially}, that is, $\B \cong \A \oplus \Cc$. By analogy with module theory, where projective and injective modules are characterized by the trivial splitting of short exact sequences, it is natural to ask whether certain Lie algebras behave like projective or injective objects.
More specifically, we consider the following questions:
\begin{enumerate}
	\item Under what conditions on $\Cc$ does every extension~\eqref{LE0} split trivially?
    
	\item Under what conditions on $\A$ does every extension~\eqref{LE0} split trivially?
\end{enumerate}

Levi's theorem~\cite[Part~I, Chapter~VI, Theorem~4.1]{Serre} ensures that whenever $\Cc$ is semisimple, the Lie algebra homomorphism $\pi$ always admits a right inverse. However, such splitting is generally semi-trivial: extension~\eqref{LE0} only splits as a semidirect sum, not necessarily as a direct sum. When $\Cc$ is not semisimple, $\pi$ may fail to have a right inverse, and extension~\eqref{LE0} cannot split even semi-trivially. Therefore, in Definition~\ref{DSE} of Section~\ref{sec2}, we distinguish two notions of splitting:
\begin{itemize}
    \item \emph{semi-trivial splitting} (or \emph{split}), meaning that $\pi$ admits a right inverse, and $\B\cong \A \oplus_\rho \Cc$ via an action $\rho$ of $\Cc$ on $\A$.
    
    \item \emph{trivial splitting}, meaning that $\B\cong \A\oplus \Cc$.
\end{itemize}
Our first main result shows that, in contrast to the module category, the category of Lie algebras admits no nontrivial objects that behave like projective or injective modules with respect to arbitrary extensions. Specifically, we have Theorem \ref{thm1} below which is proven in Section~\ref{sec2}.

\begin{theorem}\label{thm1}
    We have the following assertions.
    \begin{enumerate}
        \item\label{thm1-1} For any nontrivial Lie algebra $\Cc$, there exists an extension~\eqref{LE0} of $\Cc$ that does not split trivially, i.e., $\B \ncong \A \oplus \Cc$. 
        \item\label{thm1-2} For any nontrivial Lie algebra $\A$, there exists an extension~\eqref{LE0} by $\A$ that does not split trivially, i.e., $\B \ncong \A \oplus \Cc$.
    \end{enumerate} 
\end{theorem}

Theorem~\ref{thm1} indecates that true projectivity or injectivity cannot exist in the category of Lie algebras.
Motivated by this, we restrict our attention to two natural and symmetric subclasses of extensions: \emph{central} and \emph{cocentral} extensions.
A central extension of a Lie algebra $\Cc$ by an abelian Lie algebra $V$ is an extension
\begin{equation}\label{CE0}
    \xymatrix
	{
		0 \ar[r] & V \ar[r]^\iota & \B \ar[r]^\pi & \Cc \ar[r] & 0
	}\tag{CE}
\end{equation}
in which $\iota (V) \subset \Z(\B)$, the center of $\B$. By Levi's theorem, all central extension~\eqref{CE0} of semisimple Lie algebra split trivially, meaning that semisimple Lie algebras possess a ``projective-type'' property with respect to central extensions. However, this is a special situation where the kernel $\ker \pi = V$ is abelian. To obtain a symmetric dual notion, we introduce the following extensions
\begin{equation}\label{CCE0}
    \xymatrix
	{
	   0 \ar[r] & \A \ar[r]^\iota & \B \ar[r]^\pi & V \ar[r] & 0
	}\tag{CCE}
\end{equation}
which we refer to in this paper as the \emph{cocentral extensions} of $V$ by $\A$.

Surprisingly, the trivial splitting of cocentral extensions~\eqref{CCE0} turns out to be closely related to \emph{complete} Lie algebras, i.e., those with trivial center and only inner derivations. Historically, complete Lie algebras were introduced by Chevalley~\cite{Che44} in 1944 and formalized by Jacobson~\cite{Jac62} in 1962. Fundamental results in Lie theory show that semisimple Lie algebras are complete, and thus attention naturally turns to the solvable case. The second main result of the paper establishes an equivalence between trivial split cocentral extensions and solvable complete Lie algebras as in Theorem~\ref{thm2} whose proof is given in Section~\ref{sec3}. Combining with Proposition~\ref{P36} on the closedness under finite direct sums of complete Lie algebras, this means that solvable complete Lie algebras have a ``injective-type'' property with respect to cocentral extensions.

\begin{theorem}\label{thm2}
   For a solvable Lie algebra $\A$, the following statements are equivalent.
 \begin{enumerate}		
    \item\label{thm2-1} $\A$ is a complete.
    \item\label{thm2-2} Every cocentral extension~\eqref{CCE0} of any abelian Lie algebra $V$ by $\A$ splits trivially.
\end{enumerate}  
\end{theorem}

Dually, we introduce \emph{cocomplete} Lie algebras,
defined as those Lie algebras $\Cc$ for which every central extension~\eqref{CE0} by an abelian Lie algebra splits trivially.
This means that cocomplete Lie algebras have an appropriate analogue of ``projective-type'' property with respect to central extensions. As a dual version of Proposition \ref{P36}, we prove Proposition~\ref{P43} on the property connecting complete duality and finite direct sums. Next, we prove in Proposition~\ref{P321} that $\Cc$ is cocomplete if and only if its second cohomology group with trivial coefficients vanishes, i.e., $H^2(\Cc,\K)=0$. In particular, semisimple Lie algebras are cocomplete. Therefore, in the category of Lie algebras, the class of semisimple Lie algebras coincides with the intersection of the classes of complete and cocomplete Lie algebras.

Although the class of cocomplete Lie algebras appears to be quite interesting, their classification seems to be a nontrivial problem. A simple observation shows that low-dimensional Lie algebras are often almost abelian, and this motivates us to investigate the class of Lie algebras that are not only almost abelian but also cocomplete. Historically, almost abelian Lie algebras were introduced by Kolman~\cite{Kol65} in the classification of upper semimodular Lie algebras over a field of characteristic zero. Such Lie algebras are only abelian, almost abelian, or the 3-dimensional indecomposable simple Lie algebras. This result was later extended by Gein~\cite{Gei76}.

Recall that an $(n+1)$-dimensional non-abelian Lie algebra $\Cc$ is \emph{almost abelian} if it admits an $n$-dimensional abelian subalgebra~\cite[Definition 1]{Ave22}. Nevertheless, if such a subalgebra exits, then it must be an ideal~\cite[Proposition 3.1]{BC12}, i.e., there exists an abelian ideal $\K^n \subset \Cc$ such that $[\Cc, \K^n] \subset \K^n$. Equivalently, $\Cc$ can be viewed as the 1-dimensional extension of $\K^n$ by a derivation $D \in \Der(\K^n) = \mathfrak{gl}_n(\K)$, namely,
\[
	\Cc = \K^n \oplus_D \K e_0, \quad e_0 \in \Cc \setminus \K^n.
\]
Within this framework, our next main result presents a necessary and sufficient condition for which an almost abelian Lie algebra is cocomplete in Theorem~\ref{thm3}.

\begin{theorem}\label{thm3}
    Assume that $\K \in \{\R, \C\}$. An $(n+1)$-dimensional almost abelian Lie algebra $\Cc = \K^n \oplus_D \K e_0$ is cocomplete if and only if the following conditions hold.
\begin{enumerate}
	\item\label{thm3-1} $D$ is an automorphism of $\K^n$, i.e., $D \in \Au(\K^n) \equiv \GL_n(\K)$. 
	\item\label{thm3-2} There is no pair of $\lambda, \mu$ of $D$ (if $\K = \C$) or of the complexification $D_\C \co D \otimes _\R \id_\C$ of $D$ (if $\K = \R$) whose sum is zero.
\end{enumerate}
\end{theorem}

Theorem~\ref{thm3} provides a full classification of almost abelian cocomplete Lie algebras since it reduces the problem to classifying invertible square matrices satisfying Condition~\ref{thm3-2}, up to \emph{proportional similarity} (see Proposition~\ref{P354}). In combination with computer algebra techniques in~\cite{NLV25}, this approach yields a classification of almost abelian cocomplete Lie algebras up to isomorphism. The proof of Theorem~\ref{thm3} as well as the resulting classification is presented in Section~\ref{sec4}.

The structure of the paper is as follows. Section~\ref{sec2} presents the concept of Lie algebra extensions and proves Theorem~\ref{thm1}. Section~\ref{sec3} introduces cocentral extensions and clarifies the characterization of complete Lie algebras through Theorem~\ref{thm2}. Section~\ref{sec4} is dedicated to cocomplete Lie algebras, where we present their characterizations and the classification of almost abelian cocomplete Lie algebras in Theorem~\ref{thm3}.

\section{Extensions of Lie algebras}\label{sec2}

In this section, we first recall extensions of Lie algebras for later use. At the end of this section, we will prove Theorem~\ref{thm1} -- the first main result of the paper.

\subsection{Lie algebra extensions}

\begin{definition}[Lie algebra extensions]
    Let $\A$, $\B$ and $\Cc$ be Lie algebras.
	\begin{enumerate}
		\item A sequence of Lie algebra homomorphisms of the form~\eqref{LE0}
		\begin{equation}
		\xymatrix
	{
 		0 \ar[r] & \A \ar[r]^\iota & \B \ar[r]^\pi & \Cc \ar[r] & 0
 	}\tag{LE}          
 		\end{equation}
		is called a \emph{short exact sequence} of Lie algebras if $\iota$ is injective, $\pi$ is surjective and $\im \iota = \ker \pi$ is an ideal in $\B$. 
		In this case, the sequence~\eqref{LE0} as well as the Lie algebra $\B$ is called a \emph{Lie algebra extension} (hereafter, an \emph{extension}) of $\Cc$ by $\A$.
		\item The following extension
		\begin{equation}\label{TE1}
        \xymatrix
        {
            0 \ar[r] & \A \, \ar@{^{(}->}[r]^{\hspace{-10pt}e} & \A \oplus \Cc \ar@{->>}[r]^{\hspace{10pt}pr} & \Cc \ar[r] & 0
        }\tag{TLE}
		\end{equation}
		is called the \emph{trivial Lie algebra extension} (hereafter, \emph{trivial extension}) of $\Cc$ by $\A$. Here, $e \colon \A \hookrightarrow \A \oplus \Cc$, $a \mapsto (a, 0)$ is the canonical embedding, and $pr \colon \A \oplus \Cc \twoheadrightarrow \Cc$, $(a, c) \mapsto c$ is the canonical projection.         
	\end{enumerate}
\end{definition}

\begin{definition}[Equivalence of extensions]\label{DEE}
	Extension~\eqref{LE0} and the following extension
	\[
	\xymatrix
	{
 		0 \ar[r] & \A \ar[r]^{\bar{\iota}} & \bar{\B} \ar[r]^{\bar{\pi}} & \Cc \ar[r] & 0
 	}
	\]
	are said to be \emph{equivalent} if there exists a Lie algebra isomorphism $f \colon \B \to \bar{\B}$ making the following diagram commute
	\[
    \xymatrix
		{
			0 \ar[r] & \A \ar[r]^\iota \ar@{=}[d] & \B \ar[r]^\pi \ar[d]_{\raisebox{1.2ex}{$\scriptstyle\cong$}}^{\raisebox{1.5ex}{$\scriptstyle f$}} & \Cc \ar[r] \ar@{=}[d] & 0 \\			
			0 \ar[r] & \A \ar[r]^{\bar{\iota}} & \bar{\B} \ar[r]^{\bar{\pi}} & \Cc \ar[r] & 0.
		}
    \]
    From now on, every extension~\eqref{LE0} that is equivalent to a trivial extension~\eqref{TE1} will also be called a \emph{trivial extension} of $\Cc$ by $\A$.
\end{definition}

\begin{definition}[Lie modules]\label{D113}
	Let $\A$ and $\B$ be two Lie algebras. If there is a Lie algebra homomorphism $\rho \colon \A \to \Der (\B)$ then $\B$ becomes an \emph{$\A$-module} by the structure $\A \times \B \to \B$ as follows  
    \[
       (a, x) \mapsto a \cdot x \co \rho(a) (x), \quad \text{for all } a \in \A \text{ and } x \in \B.
    \]
    If $\rho = 0$ then $\B$ is a \emph{trivial $\A$-module}, i.e., $a \cdot x = 0$ for all $a \in \A$ and $x \in \B$.
\end{definition}

\begin{remark}\label{R114}
For extension~\eqref{LE0}, we have the following observations. 
    \begin{enumerate}
        \item If we only consider extension~\eqref{LE0} in the category of vector spaces, then obviously $\B \cong \A \dotplus \Cc$. Without loss of generality, we can consider $\B \equiv \A \dotplus \Cc$ and each element $b \in \B$ can be written uniquely in the form $(a, c) \in \A \dotplus \Cc$. Now $\iota \colon \A \hookrightarrow \B$ and $\pi \colon \B \twoheadrightarrow \Cc$ can be viewed as the canonical embedding and projection as follows
        \[
            \begin{array}{l l l l}
              \iota(a) \co (a, 0) & \text{and} & \pi(a, c) \co c, & \text{for all } a \in \A \text{ and } c\in \Cc. 
            \end{array}
        \]
        Thus, the linear embedding $\iota$ always has a left linear inverse $r \colon \B \to \A$, and the linear projection $\pi$ always has a right linear inverse $q \colon \Cc \to \B$ with
        \[
           \begin{array}{l l l l}
              r(a, c) \co a & \text{and} & q(c) \co (0, c), & \text{for all } a \in \A \text{ and } c\in \Cc.
            \end{array}
        \]
        
        \item\label{R114-2} Now we consider extension~\eqref{LE0} in the category of Lie algebras. If $r$ is also the left inverse of $\iota$ as a Lie algebra homomorphism then $r$ induces a trivial $\Cc$-module structure $\Cc \times \A \to \A$ on $\A$ as follows
        \[
            c \cdot a \co [r(0, c), a]_\A = [0, a]_\A = 0, \quad \text{for all } c \in \Cc \text{ and } a \in \A.
        \]
        Hence,~\eqref{LE0} becomes a trivial extension of $\Cc$ by $\A$, i.e., $\B = \A \oplus \Cc$ (see~\cite{Che-Eil}).
        
        \item\label{R114-3} We continue to consider extension~\eqref{LE0} in the category of Lie algebras. If $q$ is also the right inverse of $\pi$ as a Lie algebra homomorphism then it also induces a $\Cc$-module structure on $\A$, but this structure is generally not trivial. Specifically, for all $a \in \A$, $c \in \Cc$, we have 
        \[ \Cc \times \A \to \A; \, 
            (c, a) \mapsto c \cdot a \co \iota^{-1}([q(c), \iota(a)]_\B) = \iota^{-1}([(0,c), (a,0)]_\B).
        \]
        In other words, we get the action $\rho_q \colon \Cc \to \Der(\A)$ of $\Cc$ on $\A$ as follows
        \[
            \rho_q(c)(a) \co c \cdot a = \iota^{-1}([q(c), \iota(a)]_\B); \, \forall a \in \A, \, \forall c \in \Cc
        \]
        and $\B = \A \oplus_{\rho_q} \Cc$ via the action $\rho_q$ (see~\cite{Che-Eil}). That means $\B = \A \dotplus \Cc$ and for all $(a_1, c_1), (a_2, c_2) \in \B$, their Lie bracket is defined as follows
        \[
            [(a_1, c_1), (a_2, c_2)]_{\B} := \Bigl([a_1, a_2]_{\A} + \rho_{q} (c_1)(a_2) - \rho_{q}(c_2)(a_1), \, [c_1, c_2]_{\mathcal{C}}\Bigr).
        \]
        Note that if $\iota$ maps $\A$ into the center of $\B$, i.e., $\im \iota = \ker \pi \subseteq \Z(\B)$, the action $\rho _q$ induced by $q$ will be trivial and~\eqref{LE0} also becomes a trivial extension of $\Cc$ by $\A$, i.e., $\B \equiv \A \oplus \Cc$ as the direct sum of Lie algebras.
    \end{enumerate}
\end{remark}

\begin{proposition}[{see \cite{Che-Eil}}]\label{P121}
	For an extension~\eqref{LE0}, the followings are equivalent.
	\begin{enumerate}
		\item\label{P25-1} It is a trivial extension of $\Cc$ by $\A$. 
		\item\label{P25-2} The Lie algebra monomorphism $\iota$ has a left inverse $r \colon \B \to \A$.
        
        \item\label{P25-3} The Lie algebra epimorphism $\pi$ has a right inverse $q \colon \Cc \to \B$ and the action $\rho _q$ of $\Cc$ on $\A$ induced by $q$ is trivial.      
	\end{enumerate}
\end{proposition}
Before proving Theorem~\ref{thm1}, we need to define precisely the concepts of splitting and trivial splitting.

\begin{definition}[The ``trivial split'' and ``split'' of extensions]\label{DSE}
    We say that:    
    \begin{enumerate}
        \item Extension~\eqref{LE0} \emph{splits semi-trivially} (or simply, \emph{splits}) if the Lie algebra epimorphism $\pi$ has a right inverse $q$. In particular, then $\B = \A \oplus_{\rho_q} \Cc$ via the action $\rho _q$ of $\Cc$ on $\A$ induced by $q$.
        
        \item Extension~\eqref{LE0} \emph{splits trivially} if it satisfies one (and thus all) of the equivalence conditions stated in Proposition~\ref{P121}. By the first statement, when \eqref{LE0} splits trivially, we are also entitled to say that it is a {\it trivial extension} of $\Cc$ by $\A$.
    \end{enumerate}
\end{definition}

\begin{remark}\label{R123}
    For the terms ``split trivially'' and ``split'' (semi-trivially), we have a few comments below.
    \begin{enumerate}
        \item Schottenloher~\cite[Definition 4.2]{Sch08} uses the term ``split'' which is actually ``split semi-trivially'' for every extension. But the word ``split'' for any central extension is ``split trivially'' as in Definition~\ref{DSE}.
        \item Weibel~\cite[Exercise 7.6.1]{Weibel} uses the term ``split'' which is actually ``split semi-trivially'' for every abelian Lie algebra extension as in Definition~\ref{DSE}.
        \item The terms ``trivial'' and ``split'' in Kerf et al.~\cite[Definition 18.1.6]{KBK97} correspond to ``split trivially'' and ``split semi-trivially'' in Definition~\ref{DSE}, respectively.
    \end{enumerate} 
\end{remark}

\subsection{Proof of Theorem~\ref{thm1}}
As mentioned in Section~\ref{sec1}, semisimple Lie algebras closely resemble projective modules and are closely related to semi-trivially split extensions. Specifically, we have the following well-known result~\cite[Part I, Chapter VI, Theorem 4.1]{Serre}.
\begin{proposition}[Levi's theorem]\label{P124}
	Let $\pi \colon \B \to \Cc$ be a surjective homomorphism of a Lie algebra $\B$ onto a semisimple Lie algebra $\Cc$.
	Then there exists a Lie algebra homomorphism $q \colon \Cc \to \B$ such that $\pi \circ q = \id_\Cc$.
\end{proposition}
By Levi's theorem, if $\Cc$ is a semisimple Lie algebra then every extension~\eqref{LE0} always splits semi-trivially. However, if we replace ``split semi-trivially'' by ``split trivially'', Levi's theorem is no longer valid.  That is the first result of the paper which is stated in Theorem~\ref{thm1}. Here, we give the proof of this theorem.

\begin{proof}[{\bf Proof of Theorem~\ref{thm1}}]
    We will prove each assertion of the theorem.
    \begin{enumerate}
        \item Assume that $\Cc$ is a nontrivial Lie algebra. By Ado's theorem~\cite[Chapter V, \S8, p.153]{Serre}, there exists a faithful representation $\rho \colon \Cc \to \mathfrak{gl}(V) \equiv \Der(V)$ of $\Cc$ in a certain vector space $V$. Choose $\A \co V$ as an abelian Lie algebra. Since $\Cc \neq \{0\}$, the representation $\rho$ is nontrivial, i.e., there exist $a_0 \in \A$ and $c_0 \in \Cc$ such that $\rho(c_0)(a_0) \neq 0$. Then, $\rho$ defines a $\Cc$-module structure on $\A$ as follows
    \[
        \Cc \times \A \to \A, \; (c, a) \mapsto c \cdot a \co \rho (c)(a). 
    \]
    This structure is not trivial since $c_0 \cdot a_0 \neq 0$. Now we define the Lie algebra $\B \co \A \oplus_\rho \Cc$, i.e., $\B = A \dotplus \Cc$ with Lie bracket
    \[
        [(a_1,c_1),(a_2,c_2)]_\B \co (c_1 \cdot a_2 - c_2 \cdot a_1, [c_1,c_2]_\Cc)
    \]
    for $a_1, a_2 \in \A$ and $c_1, c_2 \in \Cc$. Then we obtain an extension of the form~\eqref{LE0} with $\iota \colon \A \to \B$ and $\pi \colon \B \to \Cc$ as follows
    \[
        \begin{array}{l l l}
            \iota(a) \co (a, 0) & \text{and} & \pi(a, c) \co c.
        \end{array}
    \]
    Since $[(0, c_0), (a_0, 0)] = (c_0 \cdot a_0, 0) \neq 0$, we have $[\Cc, \A] \neq \{0\}$ which implies $\B = \A \oplus_\rho \Cc \ncong \A \oplus \Cc$. Therefore, the extension~\eqref{LE0} does not split trivially.
        \item We will consider two cases where $\A$ is abelian and non-ablelian.
        \begin{itemize}
            \item {\bf Case 1: $\boldsymbol{\{0\} \ne}\boldsymbol{\mathcal{A}} \equiv V$ is abelian}. In particular, $\ad (V) = \{0\}$. Now, we choose a nonzero derivation $D \in \Der(V) \equiv \mathfrak{gl}(V)$, i.e., $D$ is an outer derivation of $V$. Let $\Cc \co \K$ and $\B \co \A \oplus_D \K$ by $D$. Obviously, we have an extension of $\K$ by $\A$ as follows
            \begin{equation}\label{S21}
                 \xymatrix
	{
 		0 \ar[r] & \A \ar[r]^e & \B \ar[r]^{pr} & \K \ar[r] & 0,
 	}
            \end{equation}
            where $e$ is the canonical embedding and $pr$ is the canonical projection. Since $D\neq 0$, the action of $\K$ on $\A$ is nontrivial. Therefore, $\B \ncong \A \oplus\K$ and extension~\eqref{S21} does not split trivially.
           
            \item {\bf Case 2: $\boldsymbol{\mathcal{A}}$ is non-abelian}. Now we choose $\Cc \co \A$ and the action of $\Cc$ on $\A$ is the adjoint representation $\ad \colon \Cc \to \Der(\A)$, $c \mapsto \ad(c) \co \ad_c$. Let $\B \co \A \oplus_\ad \Cc$. Obviously, we have an extension of $\Cc$ by $\A$ as follows
            \begin{equation}\label{S22}
             \xymatrix
	{
 		0 \ar[r] & \A \ar[r]^\iota & \B \ar[r]^\pi & \Cc \ar[r] & 0,
 	}
            \end{equation}
            where $\iota$ and $\pi$ are the canonical embedding and projection, respectively. Since $\Cc = \A$ is nontrivial, so is $\ad$. Therefore, $\B \ncong \A\oplus \Cc$ and extension~\eqref{S22} does not split trivially.
        \end{itemize}
    \end{enumerate}
    The proof of Theorem~\ref{thm1} is complete.
\end{proof}  

\section{Central - cocentral Lie algebra extensions and complete Lie algebras}\label{sec3}

Throughout the rest of this paper, $\K$ will be $\R$ or $\C$, and $\A$, $\B$, $\Cc$ are Lie algebras over $\K$. We always consider a $\K$-vector space $V$ as an abelian Lie algebra.

\subsection{Central and cocentral Lie algebra extensions}

In this subsection, we first recall the concept of a central extension. Next, we will introduce the concept of {\it cocentral extension} that arises when the central extensions are viewed from a ``symmetric'' perspective. 

\begin{definition}[{Central Lie algebra extensions~\cite[Definition 4.1]{Sch08}}]\label{DCLE}
A \emph{central Lie algebra extension} (for short, a \emph{central extension}) of $\Cc$ by $V$ is an extension
    \begin{equation}\label{CE1}
    \xymatrix
	{
 		0 \ar[r] & V \ar[r]^\iota & \B \ar[r]^\pi & \Cc \ar[r] & 0,
 	}\tag{CE}
\end{equation}
in which $\iota (V) \co \im \iota = \ker \pi$ is in the center $\Z(\B)$ of $\B$.
\end{definition}

\begin{definition}[Cocentral Lie algebra extensions] 
	We call the following extension
	\begin{equation}\label{CCE1}
	\xymatrix
	{
 		0 \ar[r] & \A \ar[r]^\iota & \B \ar[r]^\pi & V \ar[r] & 0,
 	}\tag{CCE}
	\end{equation}
	the \emph{cocentral Lie algebra extension} (for short, a \emph{cocentral extension}) of $V$ by $\A$. 
\end{definition}

\begin{remark}\label{R213}
	By Proposition~\ref{P121} and Definition~\ref{DSE}, the cocentral extension~\eqref{CCE0} \emph{splits trivially} if $\pi$ has a Lie algebra right inverse $q \colon V \to \B$ and $V$-module structure induced by $q$ on $\A$ is trivial.	
\end{remark}	

\subsection{Complete Lie algebras}
For a Lie algebra $\A$, we denote by $\ad(\A)$ the space of its inner derivations. Hence, $\Der(\A) \setminus\ad(\A)$ is precisely the set of all outer derivations of $\A$.

\begin{definition}[{Complete Lie algebras~\cite[Chapter I]{Jac62}}] 
	A Lie algebra $\A$ is called \emph{complete} if $\Z(\A) = 0$ and $\Der(\A) = \ad (\A)$.
\end{definition}

\begin{remark}\label{R35} The class of complete Lie algebras has the following obvious properties.
    \begin{enumerate}
        \item The trivial Lie algebra is complete.
        \item Every nonzero complete Lie algebra is non-abelian.
        \item All semisimple Lie algebras are always complete.
        \item All nonzero nilpotent Lie algebras are not complete.
    \end{enumerate}    
\end{remark}

\subsection{Characteristics of complete Lie algebras -- The proof of Theorem~\ref{thm2}}
In this subsection, we will prove the second main result of the paper -- Theorem~\ref{thm2} on the characterization of complete Lie algebras associated with the trivial splitting of cocentral extensions. First, we have the following statement about the stability of the class of complete Lie algebras with respect to the direct sums.

\begin{proposition}[The closedness under finite direct sums]\label{P36}
    The direct sum of two, and hence of a finite family of complete Lie algebras, is also complete.
\end{proposition}

\begin{proof}
    Let $\A \co \A_1 \oplus \A_2$ be the Lie algebra direct sum of two complete Lie algebras $\A_1$ and $\A_2$. Obviously, $\Z(\A) = \Z(\A_1) \oplus \Z(\A_2) = \{0\}$. Furthermore, if $D \in \Der(\A)$ then there exist $a = (a_1, a_2) \in \A$ with $a_1 \in \A_1$ and  $a_2 \in \A_2$ such that $D|_{\A_1} = \ad_{a_1}$ and $D|_{\A_2} = \ad_{a_2}$. Therefore, we have
    \[
        D(x) = \bigl(\ad_{a_1}(x_1), \ad_{a_2}(x_2)\bigr) = \ad_{a}(x), \quad \text{for all } x = (x_1, x_2) \in \A_1 \oplus \A_2 = \A.
    \]
    In other words, $D \in \ad(\A)$, i.e., $\Der(\A) = \ad(\A)$. So $\A$ is complete. By finite induction, we obtain the conclusion of Proposition~\ref{P36}.   
\end{proof}

Now we present the proof of Theorem~\ref{thm2}.

\begin{proof}[{\bf Proof of Theorem~\ref{thm2}}]
We will prove that (\ref{thm2-1} $\Rightarrow$ \ref{thm2-2}) and (\ref{thm2-2} $\Rightarrow$ \ref{thm2-1}).

\noindent {\bf ($\boldsymbol{\ref{thm2-1} \Rightarrow \ref{thm2-2}}$)}
Assume that $\A$ is a complete Lie algebra, i.e., its center $\Z(\A) \equiv 0$ and $\Der(\A) \equiv \ad(\A)$. We need to prove that any cocentral Lie extension~\eqref{CCE0} splits trivially. Without loss of generality, we can always consider $V \equiv \K^n$ as an $n$-dimensional vector space. In this manner, the cocentral extension~\eqref{CCE0} can be rewritten in the following form
	\begin{equation}\label{SCCE}
     \xymatrix{
	   0 \ar[r] & \A \ar[r]^\iota & \B \ar[r]^\pi & \K^n \ar[r] & 0.
	}
	\end{equation} 

By Statement 3 in Proposition~\ref{P121}, to prove that extension~\eqref{SCCE} splits trivially, we only need to show that there exists a Lie algebra homomorphism $q \colon \K^n \to \B$ which is the right inverse of $\pi$ and $q$ induces the trivial action of $\K^n$ on $\A$. Since the sequence~\eqref{SCCE} is exact, we can always view $\B \equiv \A \dotplus \K^n$ as a direct sum of vector spaces and for all $a \in \A$, $(a, v) \in \B$, $v \in \K^n$ we have
\[\A \equiv \im \iota = \A \times \{0\}, \, \iota(a) := (a, 0); \, \{0\} \times \K^n \equiv \K^n, \, \pi(a, v) := v.\]

Let \(\K^n = \s \{e_1, \ldots, e_n\} \). 
Since $\Der(\A) = \ad(\A)$, so $\Der\bigl(\A \times \{0\}\bigr) = \ad\bigl(\A \times \{0\}\bigr)$ 
and there always exists \( a_k \in \A \) such that $\ad_{(0, e_k)}|_{\A \times \{0\}} = \ad_{(a_k, 0)}$ 
for any $k = 1, \ldots, n$. 
Let $q: \K^n \to \B$ be the linear operator defined as follows: 
\[q: \K^n \to \B, \, e_k \mapsto q(e_k) = b_k := (-a_k, e_k) = (0, e_k) - (a_k, 0); \forall k = 1, \ldots n.\]

We will show that $q$ is a Lie algebra homomorphism. First, for any $k \in \{1, \ldots, n\}$ and every \( a\in \A \), we have 
\[[b_k, (a, 0)]_{\B} = [(0, e_k) - (a_k, 0), (a, 0)]_{\B} = \ad_{(0, e_k)}(a, 0) - \ad_{(a_k, 0)}(a, 0) = 0.\]
Apply the Jacobi identity for $b_k, b_h \in \B$ ($1 \leq k < h \leq n$) and $a \in \A$, we have
		 \[
		 \Bigl([[b_k, b_h], (a, 0)] + [[(a, 0), b_k], b_h] + [[b_h, (a, 0)], b_k] \equiv 0\Bigr) \Rightarrow \Bigl([[b_k, b_h], (a, 0)] = 0\Bigr).
		 \]
Therefore $[b_k, b_h] \in \Z\bigl(\A \times \{0\}\bigr) \equiv \Z(\A) = \{0\}; \, 1\leq k < h \leq n$. In other words, 
$[b_k, b_h] = [q(e_k), q(e_j)] = 0 = q([e_k, e_h])$ and $q$ is a Lie algebra homomorphism. It is easy to see $\pi \circ q (e_k) = \pi(-a_k, e_k) = e_k = \id_{\K^n}(e_k); \, k = 1, \ldots, n$ and $q$ is the right inverse of $\pi$.

On the other hand, since $[b_k, (a, 0)] = 0$ ($k = 1, \ldots, n$) so $[\K^n, \A] = 0$ and the action $\rho _q$ of $\K^n$ on $\A$ induced by $q$ is trivial. 
Hence $\B \cong \A \oplus \K^n$ and extension~\eqref{SCCE} splits trivially.

\smallskip

\noindent {\bf ($\boldsymbol{\ref{thm2-2} \Rightarrow \ref{thm2-1}}$)}
Suppose every Lie extension~\eqref{SCCE} splits trivially. We need show $\A$ is complete, i.e. $\Z(\A) = \{0\}$ and \( \mathrm{Der}(\A) = \mathrm{ad}(\A) \). That is, the proof will be divided into two parts: proving \( \mathrm{Der}(\A) = \mathrm{ad}(\A) \) and proving $\Z(\A) = \{0\}$.

\begin{itemize}
    \item {\bf Part 1. Show all derivations of $\A$ are inner: $\boldsymbol{\Der(\A) = \ad(\A)}$}.
    
	Suppose the opposite that $\Der(\A) \neq \ad(\A)$. Let \( D \in \mathrm{Der}(\A) \setminus \mathrm{ad}(\A) \). Define \( \B = \A \oplus _D \mathbb{K} \). In other words, \( \B = \A \dotplus  \mathbb{K} \) as a direct sum of two vector spaces and represents each element $b$ of $B$ as $b = (x, \alpha)$ with $x \in \A$ and $\alpha \in \K$. Now, the  Lie bracket on $\B$ is defined as follows
	\[[(x, \alpha), (y, \beta)] _{\B} := \bigl({[x, y]}_{\A} + D(\alpha y - \beta x), 0\bigr) \in \B; \, \forall (x, \alpha), (y, \beta) \in \B.\]
Once again, we choose the natural embedding $e$ and projection $\pi$ as follows
	\[
        e \colon \A \hookrightarrow \B, \, x \mapsto e(x) \co (x, 0); \, pr \colon \B \twoheadrightarrow \K, (x, \alpha) \mapsto pr (x, \alpha) \co \alpha.
    \]	
	Then, obviously we have
	\begin{equation}\label{NSE}
        \xymatrix
	{
 		0 \ar[r] & \A \ar[r]^e & \B \ar[r]^{pr} & \K \ar[r] & 0
 	}
	\end{equation}
	is an exact sequence. Since $D \in \Der(\A) \setminus \ad(\A)$ is an outer derivation and $\A$ is solvable, by \cite[Lemma 2.1]{Gra05} and \cite[Proposition 2.6]{Le},
    $\B = \A \oplus _D \K \ncong \A \oplus \K$, i.e., the Lie extension~\eqref{NSE} does not split trivially. This contradicts the hypothesis. Thus, $\Der(\A) = \ad(\A)$.
    
	\item {\bf Part 2. Show $\A$ has the trivial center: $\boldsymbol{\Z(\A) = \{0\}}$}.

    Suppose the opposite that $\Z(\A) \neq \{0\}$ and $z_0 \in \Z(\A) \setminus \{0\}$. In particular, $\A \ne \{0\}$. Note that if $\{0\} \ne \A \equiv V$ is commutative then $\Der(\A) \equiv \mathfrak{gl}(V) \neq \{0\} \equiv \ad(V)$. Hence, the equality $\Der(\A) = \ad(\A)$ in Part 1 implies that $\A$ is non-abelian. Now we will point-out the existence of an exact sequence of the form~\eqref{CCE0} that does not split trivially.
        
Since $\A$ is solvable and non-abelian, so $\{0\} \ne \A^1 := [\A, \A] \subsetneq \A$, in particular $\A/\A^1\neq 0$. Choose a nontrivial linear form $\omega \in (\A/\A^1)^*$ and set $\rho \co \omega \circ p \in \A^*$, where $p: \A\to \A/\A^1$ is the canonical projection. Then $\rho \ne 0$ and $\rho|_{\A^1} \equiv 0$. Now we consider the map as follows
\[D: \A \to \A; x \mapsto D(x):= \rho (x)z_0.\]
It is clear that $D$ is linear operator on $\A$ and $D|_{\A^1} \equiv 0$. On the other hand, since $z_0 \in Z(\A)$ so we also have
\[[D(x), y] = [\rho (x)z_0, y] = 0; \, \, [x, D(y)] = [x, \rho(y)z_0] = 0; \, \forall x, y \in \A.\]
Thus
\[D([x, y]) \equiv  0 \equiv [D(x), y] + [x, D(y)]; \, \, \forall x, y \in \A\]         
and $D \in \Der(\A)$. Now, let us consider $\B = \A \oplus _D \K$ is the semidirect sum of $\A$ and $\K$ 
via the derivation $D$. That means $\B = \A \dotplus \K$ as direct sum of two vector spaces $\A$ and $\K$ and the Lie bracket on $\B$ is defined as follows
        \[[(x, \alpha), (y, \beta)]_{\B} := {[x, y]}_{\A} + (\alpha \rho (y) - \beta \rho (x) ) z_0; \, \, \forall x, y \in \A; \, \, \alpha, \beta \in \K.\]
Now, we choose the natural embedding $e$ and projection $pr$ as follows 
    	\[e \colon: \A \hookrightarrow \B, \, x \mapsto e(x) \co (x, 0); \, pr \colon \B \twoheadrightarrow \K, (x, \alpha) \mapsto pr (x, \alpha) \co \alpha.\]
    Then we get the following extension
    \begin{equation}\label{S33}
        \xymatrix
	{
 		0 \ar[r] & \A \ar[r]^e & \B \ar[r]^{pr} & \K \ar[r] & 0.
 	}
    \end{equation}
        Since $\rho \ne 0$, $[\A, \K]$ must be different from 0, i.e., $\B \ncong \A \oplus \K$. Therefore, extension~\eqref{S33} does not split trivially that contradicts the hypothesis. This implies \( \Z(\A) = \{0\} \).
\end{itemize}
The proof of Theorem~\ref{thm2} is complete.
\end{proof}

\begin{remark}\label{R35}
We have some observations as follows.
\begin{itemize}
    \item As mentioned in Section~\ref{sec1}, although in the category of Lie algebras there is no object which is similar to injective objects in the category of modules, complete Lie algebras can be called ``quasi-injective'' objects, i.e., complete Lie algebras have the ``injective-type'' property with respect to cocentral extensions.
    
    \item In the proof of Theorem~\ref{thm2} above, Step {\bf ($\boldsymbol{\ref{thm2-1} \Rightarrow \ref{thm2-2}}$)} does not use the assumption on the solvability of $\A$, however, Step {\bf ($\boldsymbol{\ref{thm2-2} \Rightarrow \ref{thm2-1}}$)} does. Therefore, it is essentially possible to restate a new version of Theorem~\ref{thm2} as below.
\end{itemize}

\noindent {\bf An alternative of Theorem~\ref{thm2}.
    {\it For any Lie algebra $\A$, the following statements are equivalent.
 \begin{enumerate}		
    \item If $\A$ is a complete, then every cocentral extension~\eqref{CCE0} of any abelian Lie algebra $V$ by $\A$ splits trivially. 
    \item If $\A$ is solvable and every cocentral extension~\eqref{CCE0} of any abelian Lie algebra $V$ by $\A$ splits trivially, then $\A$ is a complete. 
\end{enumerate}}}
\end{remark}

Let $\ad \colon \A \to \mathfrak{gl}(\A)$ be the adjoint representation of $\A$. Then $\A$ becomes an $\A$-module under the adjoint action $\A \times \A \to \A$, $(x, y) \mapsto x \cdot y \co \ad_x(y) = [x, y]$. Under simple computations, the formula of coboundary operator in~\cite[Chapter IV, Section 23, formula (23.1)]{Che-Eil} yeilds 
\[
	\begin{array}{l l l}
		H^0 (\A,\A) = \Z(\A) & \text{and} & H^1(\A, \A) = \Der (\A) / \ad(\A).
	\end{array}
\]
From this point of view, complete Lie algebras can be defined equivalently by cohomologies as follows.

\begin{corollary}\label{prop-H0andH1}
	A Lie algebra $\A$ is complete if and only if its first two cohomology groups with coefficients in itself vanish, i.e., $H^0(\A, \A) = H^1(\A, \A) = 0$.
\end{corollary}

Any nonzero nilpotent Lie algebra is not complete since it has always nonzero centers. To the best of our knowledge, the classification of complete Lie algebras remains
an open problem. Here we construct an algorithm to determine whether an $n$-dimensional Lie algebra $\A = \s \{e_1, \ldots, e_n\}$ with structure constants $c_{ij}^k$ (where $1 \leq i < j \leq n$ and $k = 1, \ldots, n$) is complete. 

For $x =\sum \limits_{i=1}^n \alpha_i e_i \in \A$, we have $[x, e_j] = \sum \limits_{i,k=1}^n \alpha_i c_{ij}^ke_k$ with $j = 1, \ldots, n$. 
Therefore, $x \in \Z(\A)$ if and only if $(\alpha_1, \ldots, \alpha_n)$ is a solution of the linear system $\sum \limits_{i,k=1}^n c_{ij}^k \alpha_i  = 0$ ($j = 1, \ldots, n$). This means that $\Z(\A) = \{0\}$ if and only if this linear system has only trivial solutions. Equivalently, 
\begin{equation}\label{eq-center}
	\det \left(\sum \limits_{k=1}^n c_{ij}^k\right)_{1 \leq i,j \leq n} \neq 0 
\end{equation}

Now, assume that \eqref{eq-center}~holds. Since $\left\lbrace \ad_{e_1}, \ldots, \ad_{e_n} \right\rbrace$ is the generating set
of $\ad(\A)$, we must have $\dim \ad(\A) = n$ to ensure that $\Z(\A) = 0$. Therefore, $\A$ is complete if and only if $\Der(\A) = \ad(\A)$, i.e.,
\begin{equation}\label{der=ad}
	\dim \Der(\A) = n.
\end{equation}
Furthermore, for any derivation $D \in \Der(\A)$, if its  matrix with respect to the basis $(e_1, \ldots, e_n)$ is $[D] \in \M_n(\K)$ then $\dim \Der(\A)$ is equal to the number of independent parameters of $[D]$. In summary, we have Algorithm~\ref{alg1} below.

\begin{algorithm}[!h]
	\KwIn{Structure constants $c_{ij}^k \in \K$ of a Lie algebra $\A$}
	\KwOut{True if $\A$ is complete, otherwise False}
		\eIf{inequation~\eqref{eq-center}~does not hold}
			{return False}
			{
				Compute $\dim \Der(\A)$\;
				\eIf{equation~\eqref{der=ad}~does not hold}{return False}{return True}
			}
	\caption{Identification of complete Lie algebras}\label{alg1}
\end{algorithm}

Below we give an example to demonstrate Algorithm~\ref{alg1}. For a given Lie algebra, we only give nonzero Lie brackets, i.e., unspecified ones are zeros.
\begin{example}\label{Ex-aff(2)+aff(2)}
	Let
	\[
		\A = \s \{e_1, e_2, e_3, e_4 \colon [e_1, e_3] = e_1, [e_1, e_4] = -e_2, [e_2, e_3] = e_2, [e_2, e_4] = e_1\}.
	\]
	Here, $c_{13}^1 = -c_{14}^2 = c_{23}^2= c_{24}^1 = 1$. In this case, inequation~\eqref{eq-center} holds as
	\[
		\det \left(\sum \limits_{k=1}^4 c_{ij}^k\right)_{1 \leq i,j \leq 4} =
		\begin{vmatrix}
			0 & 0 & -1 & 1 \\
			0 & 0 & -1 & -1 \\
			1 & 1 & 0 & 0 \\
			-1 & 1 & 0 & 0 \\
		\end{vmatrix} = 4 \neq 0.
	\]
	Assume that $D \in \Der(\A)$ with matrix $[D] \in \M_4(\K)$. Then checking the Leibniz rule $D([e_i, e_j]) = [D(e_i), e_j] + [e_i, D(e_j)]$ 
	for $1 \leq i < j \leq 4$ yeilds
	\[
		[D] =
		\begin{bmatrix}
			a & b & c & d \\
			-b & a & d & -c \\
			0 & 0 & 0 & 0 \\
			0 & 0 & 0 & 0
		\end{bmatrix}; \, \,  a, b, c, d \in \K.
	\]
	Since $[D]$ contains four independent parameters, $\dim \Der(\A) = 4$. Therefore, $\A$ is complete as equation~\eqref{der=ad}~holds.
\end{example}

We have used Algorithm~\ref{alg1} to determine all complex and real complete Lie algebra up to dimension 4.
This result is given in Table~\ref{tab-completeLA}.

\begin{table}[!h]
	\centering
	\caption{Complex and real complete Lie algebras of dimension $\leq 4$}\label{tab-completeLA}
	\begin{tabular}{c l c}
		\hline Lie algebras & Non-zero Lie brackets & References \\
		\hline 
			$\aff$ & $[e_1, e_2] = e_2$ \\ 
		\hline 
			$\mathfrak{sl}_2(\K)$ & $[e_1, e_2] = e_2$, $[e_1, e_3] = -e_3$, $[e_2, e_3] = e_1$ \\ 
			$\mathfrak{so}_3(\R)$ & $[e_1, e_2] = e_3$, $[e_1, e_3] = -e_2$, $[e_2, e_3] = e_1$ \\ 
		\hline
			$\aff \oplus \aff$ & $[e_1, e_2] = e_2, [e_3, e_4] = e_4$ \\ 
			$\Sf_{4,12}$ & $[e_1, e_3] = e_1, [e_1, e_4] = -e_2, [e_2, e_3] = e_2, [e_2, e_4] = e_1$ & \cite[Sec.~17.4]{SW14} \\ 
			$\mathfrak{sl}_2(\K) \oplus \K$ \\
			$\mathfrak{so}_3(\R) \oplus \R$ \\
		\hline 
	\end{tabular}
\end{table}

\section{Cocomplete Lie algebras}\label{sec4}
In Section~\ref{sec3} we see that complete Lie algebras plays the role of ``quasi-injective'' objects in the category of Lie algebras. Section~\ref{sec4} is devoted to the introduction and study of the class of Lie algebras that plays the role of ``quasi-projective'' objects in the category of Lie algebras. That is the class of {\it cocomplete} Lie algebras. In essence, cocomplete Lie algebras have only the ``quasi-projective'' property, that means the ``projective'' property with respect to cocentral extensions.
First, we give a precise definition of cocomplete Lie algebras. Afterwards, we prove that semisimple Lie algebras are cocomplete, and establishes a necessary and sufficient condition for cocompleteness. The rest of this section is devoted to examining almost abelian cocomplete
Lie algebras that are not only almost abelian but also cocomplete.

\subsection{Cocomplete Lie algebras}

\begin{definition}[Cocomplete Lie algebras]\label{D311}
    A Lie algebra $\Cc$ is \emph{cocomplete} if every central extension~\eqref{CE1} of $\Cc$ by an abelian Lie algebra $V$ always splits trivially. 
\end{definition}

\begin{proposition}\label{P312}
   All semisimple Lie algebras are always cocomplete.
\end{proposition}

\begin{proof}
    Let $\Cc$ be a semisimple Lie algebra and \eqref{CE0} is an arbitrary central extension of $\Cc$ by $V$. By Proposition~\ref{P124}, every extension~\eqref{CE0} splits semi-trivially, i.e., the Lie algebra epimorphism $\pi$ has a right inverse $q \colon \Cc \to \B$. Furthermore, $\B = V \oplus_{\rho_q} \Cc$ with $\rho_q \colon \Cc \to \Der (\B)$ is the action of $\Cc$ on $\B$ induced by $q$. Since $V$ is abelian and $\iota (V) = \im \iota \subset \Z(\B)$, so $\rho_q$ is trivial by Item~\ref{R114-3} of Remark~\ref{R114}. Therefore, the central extension~\eqref{CE0} splits trivially, i.e., $\Cc$ is cocomplete.
\end{proof}

As a ``dual'' version of Proposition~\ref{P36}, we have the following Proposition~\ref{P43}.

\begin{proposition}\label{P43}
    If the direct sum of two Lie algebras, and thus of any finite family of Lie algebras, is cocomplete, then each summand is cocomplete.
\end{proposition}

\begin{proof}
    First, let $\Cc_1, \Cc_2$ be two Lie algebras and $\Cc = \Cc _1 \oplus \Cc _2$ is cocomplete. We need to prove that both of them are also cocomplete. Suppose by contradiction that at least one of them, say $\Cc _1$, is not cocomplete. By the counterfactual hypothesis, $\Cc _1$ admits at least one central extension that does not split trivially as follows
\begin{equation}\label{S41}
    \xymatrix
	{
 		0 \ar[r] & V_1 \ar[r]^{\iota _1} & \B _1 \ar[r]^{\pi _1} & \Cc _1 \ar[r] & 0
 	}
\end{equation}
here $\B _1 \ncong V_1 \oplus \Cc _1$ as the direct sum of Lie algebras. Now, we set $\B \co \B _1 \oplus \Cc_2$ as the direct sums of Lie algebras. 
Next, let us consider the following Lie algebra homomorphisms
\[
    \iota \co \iota _1 \oplus 0 \colon V_1 \equiv V_1 \oplus \{0\} \to \B _1 \oplus \Cc_2 = \B, \, v_1 \mapsto \iota (v_1) \co \bigl(\iota _1(v_1), 0\bigr); 
\]
\[\pi \co \pi _1 \oplus \id_{\Cc_2} \colon \B = \B _1 \oplus \Cc_2 \to \Cc _1 \oplus \Cc_2 = \Cc, \, (b_1, x) \mapsto \pi (b_1, b_2)) \co \bigl(\pi _1(b_1), b_2); \, \]
It is obvious that $\iota \co \iota _1 \oplus 0$ is monomorphic, $\pi \co \pi_1 \oplus \id_{\Cc_2}$ is epimorphic and $\im \iota = \im (\iota _1 \oplus 0) \equiv \im \iota _1 = \ker \pi _1 \equiv \ker (\pi _1 \oplus \id_{\Cc_2}) = \ker \pi \subset \Z(\B_1) \oplus \Z(\Cc_2) = \Z(\B)$. Therefore, the sequence
\begin{equation}\label{S42}
    \xymatrix
	{
 		0 \ar[r] & V_1 \ar[r]^{\iota} & \B \ar[r]^{\pi} & \Cc \ar[r] & 0
 	}
\end{equation}
is a central extension of $\Cc$ by $V_1$. For the central extension~\eqref{S42}, since $\B _1 \ncong V_1 \oplus \Cc _1$, so $\B = \B _1 \oplus \Cc_2 \ncong (V_1 \oplus \Cc _1) \oplus \Cc_2 \equiv V_1 \oplus (\Cc _1 \oplus \Cc_2) = V_1 \oplus \Cc$. In other words, the central extension~\eqref{S42} cannot split trivially, which contradicts the assumption that $\Cc$ is cocomplete. Therefore, both $\Cc _1, \Cc _2$ must be cocomplete. Finally, by finite induction, we obtain the conclusion of Proposition~\ref{P43}.  
\end{proof}

As mentioned in Section~\ref{sec1}, in the category of Lie algebras, the class of semisimple Lie algebras is the intersection of the classes of complete and cocomplete ones.
The next result provides a computationally necessary and sufficient condition for verifying whether a solvable Lie algebra $\Cc$ is cocomplete. To formulate and prove this result, we first need Lemma~\ref{H^p(g,V)andH^p(g,K)} as follows.

\begin{lemma}\label{H^p(g,V)andH^p(g,K)}
	For any Lie algebra $\Cc$ and every vector space $V$ as a trivial $\Cc$-module, we have the following isomorphism of vector spaces
	\[
		H^p(\Cc, V) \cong H^p(\Cc, \K) \otimes_\K V, \quad (p \in \N).
	\]
\end{lemma}

To prove the lemma, we first recall the traditional notation (see, e.g.,~\cite{CCL99,KBK97}). For vector spaces $V$ as trivial $\Cc$-modules and $0 < p \in \N$, we denote $C^0(\Cc, V) \co V$, $\Lambda^1(V^*) \co V^*$ and 
    \[
       \begin{array}{c}
         C^p(\Cc, V) \co \left\lbrace f \colon \underbrace{\Cc \times \cdots \times \Cc}_p \to V \,\Big|\, f \text{ is alternating $p$-linear} \right\rbrace, \\
         \Lambda^p (V^*) \co C^p(V, \K) = \{\omega \,| \, \omega \text{ is a $p$-form on } V\}.
       \end{array}
    \]
    
\begin{proof}[Proof of Lemma~\ref{H^p(g,V)andH^p(g,K)}]
Since $V$ is a trivial $\Cc$-module, we have
 	\begin{equation}\label{Cp(g,V)}
 		C^p(\Cc, V) \cong C^p(\Cc, \K) \otimes_\K V.
 	\end{equation}
 	Let $d^p_V \colon C^p(\Cc, V) \to C^{p+1}(\Cc, V)$, $d^p_{\K} \colon \Lambda^p (\Cc^*) \to \Lambda^{p+1}(\Cc^*)$
    be coboundary operators. Then we have $d^p_V = d^p_{\K} \otimes \id_V$ which implies 
    \[
      \begin{cases}
        Z^p(\Cc, V) \equiv \ker d^p_V \cong \ker d^p_{\K}\otimes_\K V, \\
        B^p(\Cc, V) \equiv \im d^{p-1}_V \cong \im d^{p-1}_{\K}\otimes_\K V, & p > 0.     
      \end{cases}
    \]
Therefore, we get
 	\[
 		H^p(\Cc, V) = \frac{Z^p(\Cc, V)}{B^p(\Cc, V)} \cong \frac{\ker d^p_{\K}\otimes_\K V}{\im d^{p-1}_{\K}\otimes_\K V}
 		\cong \frac{\ker d^p_\K}{\im d^{p-1}_\K} \otimes_\K V = H^p(\Cc, \K) \otimes_\K V
 	\]
	which completes the proof.
\end{proof}

Now we prove that the cocompleteness of a Lie algebra is equivalent to the vanishing of its second cohomology group with trivial coefficients as follows.

\begin{proposition}\label{P321}%
   A Lie algebra $\Cc$ is cocomplete if and only if $H^2(\Cc, \K) = 0$.
\end{proposition}

\begin{proof}
    According to Schottenloher~\cite[Chapter IV, Remark 4.7]{Sch08}, the central extension~\eqref{CE0} splits trivially if and only if $H^2(\Cc, V) = 0$. Therefore, by Lemma~\ref{H^p(g,V)andH^p(g,K)}, $\Cc$ is cocomplete if and only if 
	\[
		H^2(\Cc, V) = H^2(\Cc, \K) \otimes_\K V = 0 \iff H^2(\Cc, \K) = 0
	\]
	which completes the proof.
\end{proof}

\begin{remark}
      By means of Proposition~\ref{P321}, we not only obtain a simpler proof of Proposition~\ref{P43} but also strengthen its result. First of all, one has the K\"unneth formula~\cite[Exercise~7.3.8]{Weibel} as follows
      \begin{equation}\label{Kunneth}
          H^2(\Cc_1 \oplus \Cc_2, \K) \cong H^2(\Cc_1, \K) \oplus \left[H^1(\Cc_1, \K) \otimes H^1(\Cc_2, \K)\right] \oplus H^2(\Cc_2, \K).
      \end{equation}
      If $\Cc_1 \oplus \Cc_2$ is cocomplete, i.e., $H^2(\Cc_1 \oplus \Cc_2, \K) = 0$, by Proposition~\ref{P321}, we have 
      \[
        H^2(\Cc_1, \K) = H^1(\Cc_1, \K) \otimes H^1(\Cc_2, \K) = H^2(\Cc_2, \K) = 0.
    \]
    Therefore, $\Cc_1$ and $\Cc_2$ are cocomplete, i.e., we regain the conclusion of Proposition~\ref{P43}. In addition, the vanishing of the second summand above means that $H^1(\Cc_1, \K) = 0$ or $H^1(\Cc_2, \K) = 0$. Since $H^1(\Cc, \K) \cong (\Cc/[\Cc, \Cc])^*$, we must have $\Cc_1 = [\Cc_1, \Cc_1]$ or $\Cc_2 = [\Cc_2, \Cc_2]$, i.e., $\Cc_1$ or $\Cc_2$ is perfect Lie algebra (semisimple Lie algebras are perfect).
\end{remark}

In Lie algebra cohomologies, it is well-known that $b_p(\Cc) \co \dim H^p(\Cc, \K)$ is called the \emph{$p^{\rm th}$ Betti number} of $\Cc$.
Therefore, we obtain the following.

\begin{corollary}
	A Lie algebra $\Cc$ is cocomplete if and only if $b_2(\Cc) = 0$. 
\end{corollary}

\begin{remark}
	According to Dixmier~\cite[Th\'{e}or\`{e}me 2]{Dix55}, if $\Cc$ is nilpotent then it is not cocomplete, since $b_2(\Cc) \geq 2$.
\end{remark}

The classification of cocomplete Lie algebras is not straightforward and will be a subject of future work. 
Here, we construct an algorithm to identify whether an $n$-dimensional Lie algebra $\Cc = \s \{e_1, \ldots, e_n\}$
with structure constants $c_{ij}^k$ is cocomplete. First of all, let $\Lambda^1 (\Cc^*) \co \Cc^* = \s \{e_1^*, \ldots, e_1^*\}$ be the dual space of $\Cc$. Then, the exterior differential (or coboundary operator) of $e_i^*$ is
\begin{equation}\label{exterior-derivative}
	de_i^* = -\frac{1}{2} \sum \limits_{j,k=1}^n c_{jk}^i e_j^* \wedge e_k^*.
\end{equation}
Next, the space of 2-forms is $\Lambda^2 (\Cc^*) = \s \left\lbrace e_i^* \wedge e_j^* \,|\, 1 \leq i < j \leq n \right\rbrace$ with
\begin{equation}\label{second-exterior-derivative}
	d(e_i^* \wedge e_j^*) = de_i^* \wedge e_j^*- e_i^* \wedge de_j^*.
\end{equation}
A 2-form $\omega \in \Lambda^2 (\Cc^*)$ is \emph{closed} if $d\omega = 0$, and \emph{exact} if $\omega = d\theta$ for some $\theta \in \Lambda^1 (\Cc^*)$. Then, we have
\begin{equation}\label{Z2(g)andB2(g)}
	\begin{cases}
		Z^2(\Cc, \K) = \left\lbrace \omega \in \Lambda^2 (\Cc^*) \,|\, \omega \text{ is closed} \right\rbrace, \\
		B^2(\Cc, \K) = \left\lbrace \omega \in \Lambda^2 (\Cc^*) \,|\, \omega \text{ is exact} \right\rbrace.
	\end{cases}
\end{equation}
Due to Proposition~\ref{P321}, $\Cc$ is cocomplete if and only if $H^2(\Cc, \K) = 0$, i.e., $Z^2(\Cc, \K) = B^2(\Cc, \K)$.
This forms the basis for Algorithm~\ref{alg2} below.

\begin{algorithm}[h]
	\KwIn{Structure constants $c_{ij}^k \in \K$ of a Lie algebra $\Cc$}
	\KwOut{True if $\Cc$ is cocomplete, otherwise False}
		Compute $Z^2(\Cc, \K)$ and $B^2(\Cc, \K)$ by~\eqref{Z2(g)andB2(g)}\;
		\eIf{$Z^2(\Cc, \K) = B^2(\Cc, \K)$}
			{return True}
			{return False}
	\caption{Identification of cocomplete Lie algebras}\label{alg2}
\end{algorithm}

\begin{example}\label{Ex-cocomplete}
	Let
	\[
		\Cc = \s \{e_1, e_2, e_3, e_4 \colon [e_1, e_4] = 2e_1, [e_2, e_4] = e_2, [e_3, e_4] = e_2 + e_3\}.
	\]
	Here, $c_{14}^1 = 2$ and $c_{24}^2 = c_{34}^2 = c_{34}^3 = 1$. In this case, $\Lambda^1 (\Cc^*) = \s \{e_1^*, e_2^*, e_3^*, e_4^*\}$.
	Applying formula~\eqref{exterior-derivative}, we have 
    \[
        de_1^* = -2 e_1^* \wedge e_4^*, \; de_2^* = -e_2^* \wedge e_4^* - e_3^* \wedge e_4^*, \; de_3^* = -e_3^* \wedge e_4^*, \; de_4^* = 0.
    \]
	For $\omega = \sum \limits_{i=1}^4 a_ie_i^* \in \Lambda^1 (\Cc^*)$, one has $d\omega = -2a_1 e_1^* \wedge e_4^* - a_2e_2^* \wedge e_4^* -(a_2+a_3) e_3^* \wedge e_4^*$, and \eqref{Z2(g)andB2(g)} implies
	\[
		B^2(\Cc, \K) = \s \{e_1^* \wedge e_4^*, e_2^* \wedge e_4^*, e_3^* \wedge e_4^*\}.
	\]
	Next, we have $\Lambda^2 (\Cc^*) = \s \{e_1^* \wedge e_2^*, e_1^* \wedge e_3^*, e_1^* \wedge e_4^*, e_2^* \wedge e_3^*, e_2^* \wedge e_4^*, e_3^* \wedge e_4^*\}$. Using formula~\eqref{second-exterior-derivative}, we obtain $d(e_1^* \wedge e_4^*) = d(e_2^* \wedge e_4^*) = d(e_3^* \wedge e_4^*) = 0$ and the others $d(e_i^* \wedge e_j^*) \neq 0$. Thus, \eqref{Z2(g)andB2(g)} implies
	\[
		Z^2(\Cc, \K) = \s \{e_1^* \wedge e_4^*, e_2^* \wedge e_4^*, e_3^* \wedge e_4^*\} = B^2(\Cc, \K).
	\]
	Therefore, we conclude that $\Cc$ is cocomplete.
\end{example}

We have used Algorithm~\ref{alg2} to determine all complex and real cocomplete Lie algebras up to dimension 4.
This result is presented in Table~\ref{tab-cocompleteLA}.

\begin{table}[h]
	\centering
	\caption{Complex and real cocomplete Lie algebras of dimension $\leq 4$}\label{tab-cocompleteLA}
	\begin{tabular}{c l c}
	\hline Lie algebras & Non-zero Lie brackets & References \\
		\hline 
			$\aff$ & $[e_1, e_2] = e_2$ \\ 
		\hline 
			$\mathfrak{sl}_2(\K)$ & $[e_1, e_2] = e_2$, $[e_1, e_3] = -e_3$, $[e_2, e_3] = e_1$ \\ 
			$\mathfrak{so}_3(\R)$ & $[e_1, e_2] = e_3$, $[e_1, e_3] = -e_2$, $[e_2, e_3] = e_1$ \\ 
			$\Sf_{3,1}$ & $[e_1, e_3] = e_1$, $[e_2, e_3] = ae_2$ & \cite[Sec.~16.4]{SW14} \\ 
			$\Sf_{3,2}$ & $[e_1, e_3] = e_1$, $[e_2, e_3] = e_1+e_2$ & \cite[Sec.~16.4]{SW14} \\ 
			$\Sf_{3,3} \,(\text{over }\R)$ & $[e_1, e_3] = ae_1 - e_2$, $[e_2, e_3] = e_1 + ae_2$ & \cite[Sec.~16.4]{SW14} \\
		\hline
			$\mathfrak{sl}_2(\K) \oplus \K$ \\
			$\mathfrak{so}_3(\R) \oplus \R$ \\
			$\Sf_{4,2}$ & $[e_1, e_4] = e_1$, $[e_2, e_4] = e_1+e_2$, $[e_3, e_4] = e_2+e_3$ & \cite[Sec.~17.2]{SW14} \\
			$\Sf_{4,3}$ & $[e_1, e_4] = e_1, [e_2, e_4] = ae_2, [e_3, e_4] = be_3$ & \cite[Sec.~17.2]{SW14} \\
			$\Sf_{4,4}$ & $[e_1, e_4] = e_1, [e_2, e_4] = e_1 + e_2, [e_3, e_4] = ae_3$ & \cite[Sec.~17.2]{SW14} \\
			$\Sf_{4,5} \,(\text{over }\R)$ & $[e_1, e_4] = ae_1, [e_2, e_4] = be_2 - e_3, [e_3, e_4] = e_2 + be_3$ & \cite[Sec.~17.2]{SW14} \\
			$\Sf_{4,8}$ & $[e_1, e_4] = (1+a) e_1, [e_2, e_4] = e_2, [e_3, e_4] = ae_3$ & \cite[Sec.~17.3]{SW14} \\
			$\Sf_{4,9} \,(\text{over }\R)$ & $[e_1, e_4] = 2ae_1, [e_2, e_4] = ae_2 - e_3, [e_3, e_4] = e_2 + ae_3$ & \cite[Sec.~17.3]{SW14} \\
			$\Sf_{4,10}$ & $[e_1, e_4] = 2e_1, [e_2, e_4] = e_2, [e_3, e_4] = e_2+e_3$ & \cite[Sec.~17.3]{SW14} \\
		\hline 
	\end{tabular}
\end{table}

\subsection{Almost abelian cocomplete Lie algebras - The proof of Theorem~\ref{thm3}}

Almost abelian Lie algebras are originally defined as those admitting a basis $(e_0, e_1, \ldots, e_n)$ such that
$[e_0, e_i] = e_i$ for $i = 1, \ldots, n$ and all other brackets vanish, i.e., $\ad_{e_0} = \id_{\K^n}$. 
However, this definition is overly restrictive, as it excludes many cases where $\ad_{e_0} \neq \id_{\K^n}$. 
In this work, we adopt the broader definition as follows.

\begin{definition}[{Almost abelian Lie algebras~\cite[Definition 1]{Ave22}}]\label{D341}
    A non-abelian Lie algebra $\Cc$ is called \emph{almost abelian} Lie algebras (or for short, an AALA) if it contains an abelian subalgebra of codimension 1. When $\Cc$ is a $n$-dimensional AALA ($1 < n \in \N$), it is also called an $n$-AALA.
\end{definition}

\begin{remark}\label{R342}
    We have two remarks as follows.
    \begin{enumerate}
        \item By~\cite[Proposition 3.1]{BC12}, the subalgebra in Definition~\ref{D341}, if any, is also an ideal of $\Cc$. Moreover, every AALA $\Cc$ is 2-step solvable, i.e., $[\Cc, \Cc]$ is abelian. 
        
        \item When $\Cc$ is an $(n+1)$-AALA, we can always consider $\K^n$ as the 1-codimensional ideal of $\Cc$. Therefore, $\Cc$ is always a nontrivial 1-dimensional extension of $\K ^n$ by some derivation  $D \in \Der(\K^n) \equiv \mathfrak{gl}_n(\K)$, i.e., $\Cc = \K^n \oplus_D \K e_0$ with $e_0 \in \Cc \setminus \K ^n$. That means $\Cc \equiv \K^n \dotplus \K e_0$ and the Lie bracket of $e_0$ with any element $v \in \K ^n$ is defined by $[e_0, v] \co D(v) \in \K ^n$.
    \end{enumerate}
\end{remark}

As mentioned at the end of Section~\ref{sec1}, we will now prove Theorem~\ref{thm3} which gives a necessary and sufficient condition for which an $(n+1)$-AALA $\Cc = \K^n \oplus_D \K e_0$ is cocomplete. To this end, we need the following lemma.

\begin{lemma}\label{L343}
    If $\Cc = \K^n \oplus_D \K e_0$ is an AALA then
    \[
        H^2(\Cc,\K) \cong \Lambda^2_{D} ((\K^n)^*) \dotplus \Lambda^1 ((\K^n)^*) / \im D^T.
    \]
    Here, $D^T in \Der((\K^n)^*)$ is the transpose of $D$, and $\Lambda^2_D ((\K^n)^*)$ denotes the subspace of 2-forms $\omega \in \Lambda^2 ((\K^n)^*)$ that are cancelled by $D$, i.e.,
	\[
		D^T \cdot \omega \equiv 0 \Leftrightarrow \left(D^T \cdot \omega \right)(v_1,v_2) \co \omega \left(D(v_1),v_2 \right) + \omega \left(v_1,D(v_2)\right) = 0, \, \forall v_1, v_2 \in \K^n.
	\]
\end{lemma}

\begin{proof} 
Note that $\K^n$ has the canonical basis $(e_1, \ldots, e_n)$. So $\Cc = \K^n \dotplus \K e_0 \equiv \K ^{n+1}$ admits the basis $B \co (e_0, e_1, \ldots, e_n)$. Let $(e_1^*, \ldots, e_n^*)$ be dual basis of $(e_1, \ldots, e_n)$ in dual space $(\K^n)^*$. Then, $\Cc^* \cong (\K^{n+1})^*$ admits the dual basis $B^* \co (e_0^*, e_1 ^*, \ldots, e_n ^*)$. For any 1-form $\eta \in \Lambda^1 ((\K^n)^*)$, the Hochschild-Serre formula~\cite[p.592]{HS53} gives the 2-coboundary of $\eta$ as follows
\[
    d\eta = -e_0^* \wedge D^T \eta \in B^2 (\Cc, \K) \subset C^2(\Cc, \K).
\]
In other words, the space of 2-coboundaries
on $\Cc$ is given by
\[B^2 (\Cc, \K) = \s \{e_0 ^* \} \wedge \im D^T.\]
Now we consider the space $C^2(\Cc, \K)$ of 2-cochains on $\Cc$. Obviously, every 2-cochain $\omega \in C^2(\Cc, \K)$ on $\Cc$ can be represented as follows 
\[
    \omega = \xi + \alpha e_0^* \wedge \eta \in \Lambda^2 \bigl((\K^n)^*\bigr) \dotplus \s \{e_0^*\} \wedge \Lambda^1((\K^n)^*), 
\]
where $\alpha \in \K$. Since $\xi \in \Lambda^2((\K^n)^*)$, by the Hochschild-Serre formula, we have
\begin{align*}
    d\xi(e_0, e_i, e_j) & = -\Bigl(\xi\bigl(D(e_i), e_j \bigr) + \xi\bigl([e_i, e_j], e_0\bigr) + \xi\bigl(-D(e_j), e_i\bigr) \Bigr) \\
    & = -\Bigl(\xi\bigl(D(e_i), e_j) + \xi\bigl(e_i, D(e_j)\bigr) \Bigr) \\
    & = -\left(D^T \cdot \xi\right)(e_i, e_j), \qquad 1 \leq i < j \leq n,
\end{align*}
i.e, $d\xi(e_0, \cdot, \cdot) = -\left(D^T \cdot \xi\right)(\cdot, \cdot)$. Since $d\omega = d\xi + \alpha d(e_0^* \wedge \eta) = d\xi + \alpha \cdot 0 = d\xi$, so
\[
   d\omega = 0 \iff d\xi = 0 \iff -D^T \cdot \xi = 0 \iff \xi \in \Lambda^2_D ((\K^n)^*),
\]
which implies that the space of 2-cocycles on $\Cc$ is given by
\[
    Z^2 (\Cc, \K) = \Lambda^2_D ((\K^n)^*) \dotplus \s \{e_0^*\} \wedge \Lambda^1 ((\K^n)^*).
\]
By the definition of the second cohomology, we have
\[
    H^2(\Cc, \K) \co Z^2 (\Cc, \K) / B^2 (\Cc, \K) = \Lambda^2_D ((\K^n)^*) \dotplus \Lambda^1 ((\K^n)^*) /\im D^T.
\]
The proof of Lemma~\ref{L343} is complete.
\end{proof}

\begin{proof}[\bf Proof of Theorem~\ref{thm3}]
  According to Proposition~\ref{P321} and Lemma~\ref{L343}, an AALA $\Cc = \K^n \oplus_D \K e_0$ is cocomplete if and only if
	\[
		\Lambda^2_D ((\K^n)^*) \dotplus \Lambda^1 ((\K^n)^*)/ \im D^T = 0 \iff
		\Lambda^2_D ((\K^n)^*) = \Lambda^1 ((\K^n)^*) / \im D^T = 0.
	\]
	Naturally, we will prove that the two summands above vanish if and only if conditions \ref{thm3-1}~and~\ref{thm3-2} of Theorem~\ref{thm3} hold.
\begin{itemize}
    \item Consider the second summand $\Lambda^1 \bigl((\K^n)^*\bigr)/\im D^T$. We have
    \[
			\Lambda^1 \bigl((\K^n)^*\bigr)/\im D^T = 0 \iff 
            \im D^T = \Lambda^1 \bigl((\K^n)^*\bigr) = (\K^n)^* \iff 
            \im D = \K^n.
		\]
		This means that $D \colon \K^n \to \K^n$ is a linear isomorphism. Hence,
		\begin{equation}\label{4-1}
			\Lambda^1 \bigl((\K^n)^*\bigr)/\im D^T = 0 \iff \bigl(\text{Item~\ref{thm3-1} of Theorem~\ref{thm3} holds}\bigr).
		\end{equation}
    
    \item Consider the first summand $\Lambda^2_D \bigl((\K^n)^*\bigr)$.
		Let $\lambda, \mu$ be a pair of complex eigenvalues of $D$ (when $\K = \C$) or $D_{\C}$ (the complexification of $D$ when $\K = \R$), respectively. If $\lambda = \mu$ then $\lambda + \mu = 2\lambda \neq 0$ since $D \in \Au (\K^n) \equiv \text{GL}_n(\K)$ and $\K$ is of characteristic 0. Now, we consider the case in which $\lambda \neq \mu$, and let $v_0, w_0 \in \C^n$ be eigenvectors of $D$ (when $\K = \C$) or of the comlexification $D_\C \co D \otimes_\R \id_\C$ of $D$ (when $\K = \R$) with respect to $\lambda, \mu$. Set $\omega \co v_0^* \wedge w_0^* \in \Lambda^2 \bigl((\C^n)^*\bigr) \equiv \Lambda^2 \bigl((\R^n)^*\bigr) \otimes _\R \C$. Then $\omega \neq 0$ since $v_0$ and $w_0$ are non-proportional, and we have
		\[
  			(D^T \cdot \omega)(v_0,w_0) = -(\lambda + \mu)\omega(v_0, w_0) \, \, \, \text{or} \, \, \, (D_\C^T \cdot \omega)(v_0, w_0) = -(\lambda + \mu)\omega (v_0, w_0).
		\]
		Thus, $\omega$ is cancelled by $D$ (when $\K = \C$) or by $D_\C$ (when $\K = \R$) if and only if $\lambda + \mu = 0$.         
        Since $\Lambda^2_{D_\C} \bigl((\C^n)^*\bigr) = \Lambda^2_D \bigl((\R^n)^*\bigr) \otimes _\R \C$, we always have
		\begin{equation}\label{4-2}
			\Bigl(\Lambda^2_D \bigl((\K^n)^*\bigr) = 0 \Bigr) \wedge \bigl(\K \in \{\R, \C\}\bigr) \iff \bigl(\text{Item~\ref{thm3-2} of Theorem~\ref{thm3} holds}\bigr)
		\end{equation}
\end{itemize}
By combining \eqref{4-1} with~\eqref{4-2}, we obtain the conclusion of Theorem~\ref{thm3}.
\end{proof}

\subsection{Classification of almost abelian cocomplete Lie algebras}

Here, we solve the problem of classifying almost abelian cocomplete Lie algebras. First of all, the classification of AALAs $\Cc = \K^n \oplus_D \K e_0$ 
is well-known since it is reduced to the classification of $D \in \mathfrak{gl}_n(\K)$. Specifically, we have the following.

\begin{proposition}[{Avetisyan~\cite[Proposition 10]{Ave22}}]\label{P351}
    Two AALAs $\Cc = \K^n \oplus_D \K e_0$ and $\Cc' = \K^n \oplus_{D'} \K e'_0$ are isomorphic if and only if there exist $\alpha \in \K \setminus \{0\}$ and $\Phi \in \Au(\K^n) \equiv \GL_n(\K)$ such that $\alpha D' = \Phi^{-1} D \Phi$.
\end{proposition}

We call such two operators $D$ and $D'$ in Proposition~\ref{P351} \emph{proportionally similar}. For convenience, we recall the definition of proportional similarity below.

\begin{definition}[{Proportional similarity~\cite[Definition 2.1]{Le}}]\label{D352} 
	Two $D, D' \in \mathfrak{gl}(\K^n)$ (resp., $M, N \in \M_n(\K)$) are \emph{proportionally similar}, denoted by $D \sim _p D'$ (resp., $M \sim _p N$), if there exists $\alpha \in \K \setminus \{0\}$ and $\Phi \in \Au(\K^n) \equiv \GL_n(\K)$ (resp., $P \in \GL_n(\K)$) such that $\alpha D' = \Phi^{-1} D \Phi$ (resp., $\alpha N = P^{-1}MP$). Then, we also say $D \sim_p D'$ (resp., $M \sim _p N$) via $\alpha$ and $\Phi$ (resp., $P$).
\end{definition}

\begin{remark}\label{R353}
    Obviously, $\sim_p$ is an equivalence relation on $\M_n (\K)$. Moreover, if $M \sim_p N$ via $\alpha \in \K \setminus \{0\}$ then the Jordan canonical forms of $\alpha M$ and $N$ coincide, i.e., the sizes of their Jordan blocks are the same while eigenvalues are scaled by $\alpha$.
	Therefore, the quotient set $\M_n (\K)/\hspace{-4pt}\sim_p$ can be determined by using Jordan canonical forms of matrices.
\end{remark}

Due to Remark~\ref{R353}, Proposition~\ref{P351} reduces the problem of classifying $(n+1)$-dimensional almost abelian Lie algebras to the problem of determining the quotient set $\M_n (\K)/\hspace{-4pt}\sim_p$. In association with Theorem~\ref{thm3}, we have the following result.

\begin{proposition}\label{classification-aac}
    Let $\Cc = \K^n \oplus_D \K e_0$ and $\Cc' = \K^n \oplus_{D'} \K e_0$ be almost abelian cocomplete Lie algebras, 
    where $D, D' \in \Au(\K^n) \equiv \GL_n (\K)$ satisfy condition~\ref{thm3-2} of Theorem~\ref{thm3}. 
    Then $\Cc \cong \Cc'$ if and only if $D \sim_p D'$.
\end{proposition}

This result presents a classification of almost abelian cocomplete Lie algebras. In fact, if we denote by $\GL_n^0(\K) \subset \GL_n(\K)$ the subset of invertible matrices $D$ satisfying 
condition~\ref{thm3-2} of Theorem~\ref{thm3}, i.e., $D$ has no pair of complex eigenvalues whose sum is zero, then Proposition~\ref{classification-aac} yeilds the following.

\begin{proposition}\label{P354}
    The classification of $(n+1)$-dimensional almost abelian cocomplete Lie algebras $\Cc = \K^n \oplus_D \K e_0$ 
	is equivalent to the classification of $D \in \GL_n(\K)$ satisfying condition~\ref{thm3-2} of Theorem~\ref{thm3}
	up to proportional similarity. Equivalently, we have to determine the quotient set $\GL_n^0(\K)/\hspace{-4pt}\sim _p$
	containing equivalence classes of $\GL_n^0(\K)$ under the relation $\sim _p$.
\end{proposition}

Proposition~\ref{P354} provides the theoretical foundation for Algorithm~\ref{alg3}, which classifies almost abelian cocomplete Lie algebras. To obtain the desired classification up to isomorphism, we employ additional algorithms from~\cite{Le}.

\begin{algorithm}[h]
	\KwIn{A positive integer number $n$}
	\KwOut{A list $\Sc$ of non-isomorphic $(n+1)$-dimensional almost abelian cocomplete Lie algebras}
	$\Sc \co \{ \; \}$\;
	$\Cc_D \co \K^n \oplus_D \K e_0$ with $D \in \GL_n(\K)$\;
	$\GL_n^0(\K) \co \left\lbrace D \in \GL_n(\K) \text{ satisfies condition~\ref{thm3-2} of Theorem~\ref{thm3}} \right\rbrace$\;
	$\mathcal{D} \co \GL_n^0(\K)/\hspace{-4pt}\sim _p$\;
	\For{$D \in \mathcal{D}$}
	{
		Isomorphism verification for $\Cc_D$ by algorithms in~\cite{NLV25}\;
		$\mathcal{I} \co \left\lbrace \text{Isomorphism classes of $\Cc_D$}\right\rbrace$\;
		$\Sc \co \Sc \cup \mathcal{I}$\;
	}
	Return $\Sc$.
	\caption{Classification of almost abelian complete Lie algebras}\label{alg3}
\end{algorithm}

\begin{example}
	Let $\Cc = \K^2 \oplus_D \K e_0$ be a 3-dimensional almost abelian cocomplete Lie algebra with $D \in \Au(\K^2) \equiv \GL_2(\K)$. By condition~\ref{thm3-2} of Theorem~\ref{thm3}, $D$ does not have two non-zero complex eigenvalues whose sum is zero. Therefore, the possible Jordan canonical forms of $D$ are as follows:
	\[
		\begin{array}{l l l}
			J_1 \co
			\begin{bmatrix}
				\lambda_1 & 0 \\ 0 & \lambda_2
			\end{bmatrix}, &
			J_2 \co
			\begin{bmatrix}
				\lambda & 1 \\ 0 & \lambda
			\end{bmatrix}, &
			J_3 \co
			\begin{bmatrix}
				\lambda & 1 \\
				-1 & \lambda
			\end{bmatrix},
		\end{array}
	\]
	where $\lambda_1\lambda_2 \neq 0$ and $\lambda_1 \neq -\lambda_2$, $\lambda\neq 0$. Note that $J_3$ does not exist over $\C$ since it coincides with $J_1$. 
	By scaling, we have
	\[
		\begin{array}{l l}
			J_1 \sim_p
			\begin{bmatrix}
				1 & 0 \\ 0 & \lambda
			\end{bmatrix} \co D_1, &
			J_2 \sim_p
			\begin{bmatrix}
				1 & 1 \\ 0 & 1
			\end{bmatrix} \co D_2,
		\end{array}
	\]
	where $\lambda = \frac{\lambda_2}{\lambda_1} \neq 0, -1$.
	Hence, we have $\GL_2^0(\R)/\hspace{-3pt}\sim_p = \{D_1, D_2, J_3\}$ and $\GL_2^0(\C)/\hspace{-3pt}\sim_p = \{D_1, D_2\}$. These quotient sets create 3-dimensional almost abelian cocomplete Lie algebras with basis $\{e_0, e_1, e_2\}$ as follows:
	\[
		\begin{array}{l l l}
			\Cc_{3.1}^\lambda: & [e_0, e_1] = e_1, [e_0, e_2] = \lambda e_2 & (\lambda \neq 0, -1), \\
			\Cc_{3.2}: & [e_0, e_1] = e_1, [e_0, e_2] = e_1 + e_2, \\
			\rc_{3.3}^\lambda: & [e_0, e_1] = \lambda e_1 - e_2, [e_0, e_2] = e_1 + \lambda e_2 & (\lambda \neq 0). \\
		\end{array}
	\]
	Here, $\rc$ means this Lie algebra only exists over $\R$, others are valid over $\R$ and $\C$. 
	Finally, by using algorithms in~\cite{NLV25}, we can see that $\Cc_1^\lambda \cong \Cc_1^\frac{1}{\lambda}$ and $\rc_3^\lambda \cong \rc_3^{-\lambda}$
	by the following isomorphisms
	\[
		\begin{array}{l l l}
			\begin{bmatrix}
				\lambda & 0 & 0 \\
				0 & 0 & 1 \\
				0 & 1 & 0
			\end{bmatrix} & \text{and} &
			\begin{bmatrix}
				-1 & 0 & 0 \\
				0 & 0 & 1 \\
				0 & 1 & 0
			\end{bmatrix}
		\end{array}
	\]
	respectively, and we use the notations $\lambda \equiv \frac{1}{\lambda}$ and $\lambda \equiv -\lambda$ to indicate these phenomena.
	In summary, the class of 3-dimensional complex and real almost abelian cocomplete Lie algebras has been classified up to isomorphism.
\end{example}

We have used Algorithm~\ref{alg3} to classify up to isomorphism complex and real almost abelian cocomplete Lie algebras of dimension up to 4 over $\K$. The obtained results are given in Table~\ref{tab-aaCLA}, in which $\equiv$ indicates that the corresponding parameters give rise to isomorphic Lie algebras, and the disappearance of $\equiv$ means that the parameter is optimal in the sense that different parameters yield non-isomorphic Lie algebras.

\begin{table}[!h]
	\centering
	\caption{Complex and real almost abelian cocomplete Lie algebras of dimension $\leq 4$}\label{tab-aaCLA}
	\begin{tabular}{c l}
	\hline Lie algebras & Non-zero Lie brackets, conditions of parameters and isomorphisms \\
	\hline 
		$\aff$ & $[e_0,e_1] = e_1$ \\ 
	\hline 
		$\Cc_{3.1}^\lambda$ & $[e_0, e_1] = e_1$, $[e_0, e_2] = \lambda e_2$ \;\big($\lambda \neq 0, -1$; $\lambda \equiv \frac{1}{\lambda}$\big) \\ 
		$\Cc_{3.2}$ & $[e_0, e_1] = e_1$, $[e_0, e_2] = e_1 + e_2$ \\ 
		$\rc_{3.3}^\lambda$ & $[e_0, e_1] = \lambda e_1 - e_2$, $[e_0, e_2] = e_1 + \lambda e_2$ \;($\lambda \neq 0$; $\lambda \equiv -\lambda$) \\ 
	\hline
		\multirow{2.5}{*}{$\Cc_{4.1}^{\lambda\mu}$} & $[e_0, e_1] = e_1$, $[e_0, e_2] = \lambda e_2$, $[e_0, e_3] = \mu e_3$ 
			\;($\lambda, \mu \neq 0, -1$; $\lambda + \mu \neq 0$) \\
			& $\left((\lambda, \mu) \equiv (\mu, \lambda) \equiv (\frac{\lambda}{\mu}, \frac{1}{\mu}) \equiv (\frac{\mu}{\lambda}, \frac{1}{\lambda}) 
			\equiv (\frac{1}{\mu}, \frac{\mu}{\lambda}) \equiv (\frac{1}{\lambda}, \frac{\mu}{\lambda})\right)$ \\ 
		$\Cc_{4.2}^\lambda$ & $[e_0, e_1] = e_1$, $[e_0, e_2] = \lambda e_2$, $[e_0, e_3] = e_2 + \lambda e_3$ \;($\lambda \neq 0, -1$) \\ 
		$\rc_{4.3}^{\lambda\mu}$ & $[e_0, e_1] = \lambda e_1$, $[e_0, e_2] = \mu e_2 - e_3$, $[e_0, e_3] = e_2 + \mu e_3$ \;$(\lambda, \mu \neq 0)$ \\  
		$\Cc_{4.4}$ & $[e_0, e_1] = e_1$, $[e_0, e_2] = e_1 + e_2$, $[e_0, e_3] = e_2 + e_3$ \\ 
	\hline 
	\end{tabular}
\end{table}

\section*{Concluding remarks}

In this paper, we have introduced two new concepts: \emph{cocentral Lie algebra extensions} and \emph{cocomplete Lie algebras}. The main results can be summarized as follows.
\begin{itemize}
    \item Theorem~\ref{thm1} shows that, within the category of Lie algebras, there is no notion that exactly correspond to the projective or injective properties in the category of modules. 
    
    \item Motivated by the symmetric counterpart of central extensions, we introduced the concept of \emph{cocentral extensions}. With this notion, Theorem~\ref{thm2} characterizes the class of complete Lie algebras as those possessing an injective-type property relative to cocentral extensions
    
    \item Next, we introduce the new concept of \emph{cocomplete Lie algebras}, which play a ``projective-type'' role with respect to central extensions. 
    Specifically, every central extension of a cocomplete Lie algebra splits trivially. We also established Proposition~\ref{P321}, which provides a necessary and sufficient condition for a $\K$-Lie algebra $\Cc$ to be cocomplete, namely, $H^2(\Cc, \K) = 0$.
 
    \item Finally, Theorem~\ref{thm3} establishes the necessary and sufficient conditions for an $(n+1)$-dimensional AALA $\Cc$ to be cocomplete. 
    In addition, Proposition~\ref{P354} reduces the classification of this subclass to that of invertible $(n\times n)$-matrices up to proportional similarity. 
\end{itemize}
Our future work will focus on the classification of complete and cocomplete Lie algebras in low dimensions. In particular, we intend to describe explicitly all finite-dimensional cocomplete Lie algebras with a maximal commutative ideal of codimension 2, which will be presented in a forthcoming paper.

\section*{Acknowledgements}

This research is funded by Vietnam National University Ho Chi Minh City (VNU-HCM) under grant number B-2026-34-01.



\addcontentsline{toc}{section}{References}
\begin{thebibliography}{99}

\bibitem{Ave22}Z. Avetisyan, The structure of almost Abelian Lie algebras, Internat. J. Math. 33 (8) (2022) 2250057 (26 pages).

\bibitem{BC12}D. Burde, M. Ceballos, Abelian ideals of maximal dimension for solvable Lie algebras, J. Lie Theory 22 (3) (2012) 741--756.

\bibitem{CCL99}S. S. Chern, W. H. Chen, K. S. Lam, Lectures on Differential Geometry, Series on University Mathematics Vol. 1, World Scientific, Singapore, 1999.

\bibitem{Che44}C. Chevalley, On groups of automorphism of Lie groups,	Proc. Natl. Acad. Sci. USA 30 (1944) 274--275.

\bibitem{Che-Eil}C. Chevalley, S. Eilenberg, Cohomology theory of Lie groups and Lie algebras, Trans. Amer. Math. Soc. 63 (1) (1948) 85--124.

\bibitem{Dix55}J. Dixmier, Cohomologie des alg\`{e}bres de Lie nilpotentes, Acta Sci. Math. (Szeged) 16 (1955) 246--250.

\bibitem{Gei76}A. G. Gein, Semimodular Lie algebras, Sib. Math. J. 17 (1976) 189--193.
	
\bibitem{Gra05}W. A. de Graaf, Classification of solvable Lie algebras, Exp. Math. 14 (1) (2005) 15--25.

\bibitem{HS53} G. Hochschild and J-P. Serre, Cohomology of Lie algebras, Ann. of Math. 57 (3) (1953) 591--603.

\bibitem{Jac62}N. Jacobson, Lie Algebras, Wiley, New York, 1962.

\bibitem{KBK97}E. A. de Kerf, G. G. A. B\"auerle, A. P. E. ten Kroode, Lie Algebras, part 2: Finite and Infinite Dimensional Lie Algebras and Applications in Physics, Studies in Mathematical Physics Vol. 7, North-Holland, The Netherlands, 1997.

\bibitem{Kol65}B. Kolman, Semi-modular Lie algebras, J. Sci. Hiroshima Univ., Ser. A-1 29 (2) (1965) 149--163.

\bibitem{Le}V. A. Le, H. T. T. Cao, H. Q. Duong, T. A. Nguyen, T. N. Vo, On the problem of classifying solvable Lie algebras having small codimensional derived algebras, Comm. Algebra 50 (9) (2022) 3775--3793.

\bibitem{NLV25}T. A. Nguyen, V. A. Le, T. N. Vo, Testing isomorphism of complex Lie algebras, Comm. Algebra, doi.org/10.1080/00927872.2025.2553026

\bibitem{Sch08}M. Schottenloher, A Mathematical Introduction to Conformal Field Theory, Lecture Notes in Physics Vol. 759, Springer-Verlag, Berlin--Heidelberg, 2008.

\bibitem{Serre}J. P. Serre, Lie Algebras and Lie Groups (1964 Letures given at Harvard University), Springer-Verlag, 

\bibitem{SW14}L. \v Snobl, P. Winternitz, Classification and Identification of Lie Algebras, CRM Monograph Series Vol. 33,
American Mathematical Society, 2014.

\bibitem{Weibel}C. A. Weibel,  An Introduction to Homological Algebra, Cambridge Studies in Advanced Mathematics 38, Cambridge University Press, UK, 1994. 



\end{thebibliography}
\end{document}